\documentclass[10pt]{article}
\usepackage{amsfonts}
\usepackage{amsmath}
\usepackage{amscd,amssymb,amsthm}
\usepackage{graphics}

\newtheorem{theorem}{\sc Theorem}[section]

\newtheorem{prop}{\sc Proposition}[section]

\newcommand{\bce}{\begin{center}}
\newcommand{\ece}{\end{center}}

\begin{document}

\title{A note on a piecewise-linear Duffing-type system }

\date{}
\author{Gheorghe Tigan\thanks{Department of Mathematics, "Politehnica" University of
Timisoara, P-ta Victoriei, Nr.2, 300006, Timisoara, Timis,
Romania, email: gheorghe.tigan@mat.upt.ro }\ , Alessandro
Astolfi\thanks{ Electrical and Electronic Engineering, Imperial
College London, South Kensington Campus, London SW7 2AZ, UK,
email: a.astolfi@imperial.ac.uk } } \maketitle

\begin{abstract}

In \cite{tone2} it was suggested that the number of limit cycles
in a piecewise-linear system could be closely related to the
number of zones, that is the number of parts of the phase plane
where the system is linear. In this note we construct an example
of a class of perturbed piecewise systems with $n$ zones such that
the first variation of the displacement function is identically
zero. Then we conjecture that the system has no limit cycles using
the second variation of the displacement function expressed for
continuous functions. This system can be seen as a feedback system
in control theory.
\end{abstract}

\section{Introduction}

In this work we consider a piecewise Duffing-type system. The
classical Duffing system deals with continuous functions while the
system addressed in this paper is a discontinuous one. This system
can be seen as a feedback system in control theory. Such
non-smooth systems appear naturally in many practical systems
because many physical phenomena presents discontinuities. For
example in control theory, systems controlled by switching belong
to this class of non-smooth systems. It is known that switching
occurs in control systems in industry such as multi-body systems,
intelligent systems, robots \cite{liber}. Neuronal firing in
biology \cite{tone1} or impacts in mechanics are other systems
which address the problem. The first works on this topic are
\cite{filipov} and \cite{andronov}. Recently results can be found
in \cite{frei1}, \cite{frei2} and \cite{libre1}. In some cases,
piecewise-linear systems offer a good approximations for nonlinear
complex systems offering a valuable tool for investigating
nonlinear phenomena. As it was pointed out in \cite{libre1}, there
is a \emph{feeling} that piecewise-linear systems can present all
the features met in nonlinear dynamics, such as homoclinic or
heteroclinic orbits, limit cycles and attractors. In fact, Chua
and collaborators discovered chaotic behavior in piecewise linear
systems \cite{mada}. We will focus in the present work on the
bifurcation phenomena related to the existence of invariant closed
curves, such as limit cycles. The number and distribution
(location) of limit cycles is one of the most important problems
in qualitative theory of dynamical systems. There are numerous
works on existence, number and distribution of limit cycles for
continuous dynamical systems, for example
\cite{Cao,Tang,Li2,tig1,tig2,tig3} but not many exist for
non-smooth systems. Investigations of the dynamics of non-smooth
systems with one, two or three lines of discontinuity are
performed in \cite{gian}, \cite{frei1}, \cite{frei2}. Few papers
deal with systems with a large number of discontinuities and not
much is known about their dynamics. A recent work, which considers
a piecewise-linear function for a Li\'{e}nard system with $n$
zones is \cite{tone2}. In that work it is suggested that \emph{the
number of limit cycles is closely related to the number of zones}.
In the present work, we construct a piecewise Duffing-type system
with $n$ zones and conjecture it has no limit cycles.

\section{A piecewise-linear Duffing-type system}

Assume $H$ is a Hamilton function and consider the following
perturbed differential system

\begin{equation}\label{nr0}
\begin{array}{c}
  \dot{x}=\\
  \dot{y}=\\
\end{array}
\begin{array}{c}
\partial H/\partial y+\varepsilon f(x,y,\varepsilon) \\
-\partial H/\partial x+\varepsilon g(x,y,\varepsilon) \\
\end{array}
\end{equation}
where $f(x,y,\varepsilon)$ and $g(x,y,\varepsilon)$ are two enough
smooth functions in $x, y$ which depend analytically on a small
parameter $\varepsilon$. For any $h$ on a real interval $(a,b),$
we suppose that the set $\{ (x,y)\in\mathbb{R};H(x,y)=h \}$
contains a closed curve $C_h$ free of critical points (a circle
for example) which depends continuously on $h$. Such a family of
closed curves $C_h$ corresponds to an annulus $A$ of periodic
solutions of the unperturbed Hamiltonian system $dH=0$, that is,
the system (\ref{nr0}) for $\varepsilon=0.$\\
Recall that, if we fix a transversal segment to the flow in
(\ref{nr0}) and parameterize it using the energy level h, then the
function
\begin{equation}\label{nr01}
d(h,\varepsilon):=P(h,\varepsilon)-h=\varepsilon
M_1(h)+\varepsilon^2M_2(h)+...+\varepsilon^kM_k(h)+O(\varepsilon^{k+1}),\
h\in(a, b),
\end{equation}
where $P(h,\varepsilon)$ is the first return map (or Poincar\'{e}
map), is the displacement function defined for small
$\varepsilon.$ The variations of the displacement function,
$M_k(h)$, are also called the Melnikov functions. Computing
explicitly the Melnikov functions is a challenging problem and as
long as we know they are determined only in some particular cases.
It is known that the number of zeros of the first non-disappearing
$k-th$ order Melnikov function $M_k(h)$ provides the upper bound
of the number of limit cycles of the perturbed system emerging
from the periodic orbits of the unperturbed system \cite{za1}.
More exactly, if the first not identically null Melnikov function
is $M_{k}(h),$ then we have the following result:

\begin{theorem}\label{teo1}
If $M_1(h)=...=M_{k-1}(h)\equiv 0,$\ $M_k(h)\neq 0$ for some $h$
and $h_{1}$ is a root of the Melnikov function $M_k(h)$ such that
the $m-$derivative $M_{k}^{(m)}(h_{1})\neq 0$, $m\geq 1,$ then for
$\varepsilon \neq 0$ sufficiently small, system (\ref{nr0}) has
one limit cycle of multiplicity $m$ in an $O(\varepsilon)$
neighborhood of $C_{h_{1}}.$ In case that $M_k(h)\neq 0$ for any
$h$, then for $\varepsilon \neq 0$ sufficiently small, the system
(\ref{nr0}) has no limit cycles in an $O(\varepsilon)$
neighborhood of $C_{h}.$
\end{theorem}
\noindent More details can be found in \cite{blow} and
\cite{tone2}. The first Melnikov function for the system
(\ref{nr0}), as it is reported in \cite{ile1}, is given by:
\begin{equation}\label{nr02}
M_1(h)=\oint_{C(h)}g(x,y,0)dx-f(x,y,0)dy.
\end{equation}
\par In \cite{ile1} is presented a
method for determining the second Melnikov function, following an
algorithm described in \cite{fra1}, for a system of type
(\ref{nr0}) with the Hamiltonian expressed in the form
$H(x,y)=\frac{1}{2}y^2-U(x),$ where $U(x)$ is a polynomial of
degree at least $2$, in the case when the first Melnikov function
$M_1(h)\equiv 0.$ In this case, the second order Melnikov function
is given by:
\begin{equation}\label{nr05}
M_2(h)=\oint_{C_h}\frac{\partial g}{\partial
\varepsilon}(x,y,0)dx-\frac{\partial f}{\partial
\varepsilon}(x,y,0)dy,
\end{equation}

\noindent provided that $$\frac{\partial f}{\partial
x}(x,y,0)+\frac{\partial g}{\partial y}(x,y,0)=0.$$

\noindent For similar systems of the form

\begin{equation}\label{nr03}
\left(%
\begin{array}{c}
  \dot{x} \\
  \dot{y} \\
\end{array}%
\right) = \left(%
\begin{array}{c}
  f_{1}(x,y) \\
  f_{2}(x,y) \\
\end{array}%
\right)+\varepsilon\left(%
\begin{array}{c}
  g_{1}(x,y,\varepsilon) \\
  g_{2}(x,y,\varepsilon) \\
\end{array}%
\right),
\end{equation}
\noindent where $(x,y)\in\mathbb{R}^2,\ 0<\varepsilon\ll 1$ and
$f=(f_{1}(x,y),f_{2}(x,y))$,
$g=(g_{1}(x,y,\varepsilon),g_{2}(x,y,\varepsilon))$ are two
sufficiently smooth functions, the first Melnikov function is
treated in \cite{blow} (used in \cite{tone2}) and is given by:
\begin{equation}\label{nr04}
M_1(r)=\int_{0}^{T_r}e^{-\int_{0}^{t}div(f(\tau_{r}(s)))ds}(f_{1}(\tau_{r}(t))g_{2}(\tau_{r}(t),0)-
g_{1}(\tau_{r}(t),0)f_{2}(\tau_{r}(t)))dt
\end{equation}
\noindent where we assumed that the system (\ref{nr03}) for
$\varepsilon =0$ possesses a family of periodic orbits $C_{r}:
\tau_{r}(t), r>0$ depending on a positive real parameter $r$ and
having the period $T_{r}.$ This formula can be used even though
$f$ is differentiable but $g$ is discontinuous in some isolated
points, as remarked in \cite{tone2}. It will be employed in the
present work in such a case.\par Consider now the planar
discontinuous Duffing-type system given by:

\begin{equation}\label{nr1}
\left(%
\begin{array}{c}
  \dot{x} \\
  \dot{y} \\
\end{array}%
\right) = \left(%
\begin{array}{cc}
  0 & 1 \\
  -1 & 0 \\
\end{array}%
\right)\left(%
\begin{array}{c}
  x \\
  y \\
\end{array}%
\right)+\varepsilon\left(%
\begin{array}{c}
  0 \\
  g(x,y,\varepsilon) \\
\end{array}%
\right),
\end{equation}

\noindent where $g$ is a real function defined on the set
$I\times\mathbb{R}\times(-\varepsilon_1,\varepsilon_1),$ with
$I=(-\infty,a_1]\cup(a_1,a_2]\cup...\cup(a_n,\infty),
a_i\in\mathbb{R},i=1,2,...,n,$ $\varepsilon_1>0,$ discontinuous in
the points $a_i$ and, for a fixed $\varepsilon,$ differentiable
with respect to $(x,y)$ on each of the $(n+1)$ strips
$(a_i,a_{i+1})$ where it is defined. If $g(x,y,\varepsilon)=x^3$
and $\varepsilon = \beta$ we meet the classical Duffing system
\cite{net1} in a particular case. In \cite{tig4} we studied a
Duffing continuous system identifying the conditions to transition
to chaos.

Our first result from this paper is stated in the following
Proposition:
\begin{prop}\label{p1}
Consider two increasing sequences of real numbers \\
$a_{0}=-\infty<0<a_{1}<...<a_{n}<a_{n+1}=+\infty$,
$-\infty<\alpha_{0}<\alpha_{1}<...<\alpha_{n}<+\infty$,
$n\in\mathbb{N},n>1$ and the non-smooth linear function with $n+1$
zones
\begin{equation}\label{nr2}
g(x,y,\varepsilon)=\alpha_{i}x+\varepsilon y \ \ if \ x\in(a_{i},
a_{i+1}], \ i=0,1,...,n
\end{equation}
with the convention that the last interval is $(a_n,a_{n+1}).$

Then, for any $0<\varepsilon\ll 1$ and any two sequences
$(a_{i})_{i=0,1,...,n+1}$ and $(\alpha_{i})_{i=0,1,...,n}$ as
above, the first Melnikov function of the system (\ref{nr1}) is
identically zero.
\end{prop}
\textbf{PROOF} It is clear that the system (\ref{nr1}) is of type
(\ref{nr03}) with $f_1=y, f_2=-x, g_1=0$ and
$g_2(x,y,\varepsilon)=g(x,y,\varepsilon)=\alpha_i x+\varepsilon
y.$ With these functions, the unperturbed system (\ref{nr1}) has a
one-parameter family of periodic solutions which are circles of
the form $C_r:x^2+y^2=r^2.$ Choosing the parametrization
$\tau_r(t):x=r\sin t,\ y=r\cos t,\ t\in[0,2\pi),$ and denoting
$g(x):=g_2(x,y,0)=\alpha_{i}x,$ then using (\ref{nr04}) we get
\begin{align}\label{nr40}
M_1(r)&=\int_{0}^{2\pi}r\cos t\ g(r\sin t)dt.
\end{align}
One can observe that the same result is recovered from the formula
(\ref{nr02}), at least formally since $g$ is discontinuous,
applied for systems of the form (\ref{nr0}) with
$H(x,y)=\frac{1}{2}x^2+\frac{1}{2}y^2,$\ $f(x,y,\varepsilon)=0$
and $g(x,y,\varepsilon)=\alpha_{i}x+\varepsilon y,$ denoting the
circle $C_r:x^2+y^2=r^2$ and $g(x):=g(x,y,0)=\alpha_{i}x.$ Indeed,
from (\ref{nr02}) we have
\begin{equation}\label{nr3}
M_1(r)=\oint_{C_r}g(x)dx,
\end{equation}
and using the same parametrization $x=r\sin t,\ y=r\cos t,\
t\in[0,2\pi),$ it leads to (\ref{nr40}).

Find now $M_1(r).$ From (\ref{nr40}) we have
\begin{align}\label{nr4}
M_1(r)&=\int_{0}^{2\pi}r\cos t\ g(r\sin t)dt\\
&=\int_{0}^{\pi/2}r\cos t\ g(r\sin t)dt+\int_{\pi/2}^{3\pi/2}r\cos
t\ g(r\sin t)dt+\int_{3\pi/2}^{2\pi}r\cos t\ g(r\sin
t)dt.\nonumber \end{align}

Compute in the following the three integrals. Let
$m\in\{0,1,...,n\}$ such that $a_{m}<r<a_{m+1}$ and
$t_{0}=0<t_{1}<t_{2}<...<t_{m}<t_{m+1}=\pi/2$ an increasing
sequence of real numbers in $[0,\pi/2]$ given by $\sin
t_{i}=\frac{a_{i}}{r}, \ i=1,2,...,m.$ It is obvious now that if
$t\in(t_{i},t_{i+1}),$ then $r\sin t\in(a_{i},a_{i+1}),
i=0,1,2,...,m.$ Consequently

\begin{align}
M^{1}(r):&= \int_{0}^{\pi/2}r\cos t\ g(r\sin t)dt
=\sum^{m}_{i=0}\int^{t_{i+1}}_{t_{i}}r^2\alpha_i\cos
t\sin t dt\\
&=\frac{r^2}{2}\sum^{m}_{i=0}\alpha_i(\sin^2t_{i+1}-\sin^2t_i)\\
&=\frac{r^2}{2}\alpha_{m}-\frac{1}{2}\sum^{m}_{i=1}(\alpha_{i}-\alpha_{i-1})a_{i}^{2}\nonumber.
\end{align}

For the second integral, consider a decreasing sequence of real
numbers $t_{m+1}=\pi/2<t_{m}<t_{m-1}<...<t_{1}<t_{0}=3\pi/2$ given
by $\sin t_{i}=\frac{a_{i}}{r}, \ i=1,2,...,m.$ Because the $sin$
function is decreasing on the interval $[\pi/2,3\pi/2 ],$ one gets
that, if $t\in(t_{i+1},t_{i}),$ then $r\sin t\in(a_{i},a_{i+1}),
i=0,1,2,...,m.$ Therefore

\begin{align}
M^{2}(r):&= \int_{\pi/2}^{3\pi/2}r\cos t\ g(r\sin t)dt
=\sum^{m}_{i=0}\int^{t_{i}}_{t_{i+1}}r^2\alpha_{i}\cos
t\sin t dt\\
&=\frac{r^2}{2}\sum^{m}_{i=0}\alpha_i(\sin^2t_{i}-\sin^2t_{i+1})\\
&=\frac{r^2}{2}(\alpha_0-\alpha_{m})+\frac{1}{2}\sum^{m}_{i=1}(\alpha_{i}-\alpha_{i-1})a_{i}^{2}\nonumber.
\end{align}

In the last case, because $r\sin t<0$ for $t\in[3\pi/2,2\pi)$ we
have that $r\sin t\in(a_0,a_1)$ so that $g(r\sin t)=\alpha_0r\sin
t$.\\ Therefore

\begin{align}
M^{3}(r):= \int_{3\pi/2}^{2\pi}r\cos t\ g(r\sin t)dt=
\int_{3\pi/2}^{2\pi}r^2\alpha_0\cos t\sin t
dt=-\frac{r^2}{2}\alpha_{0} .
\end{align}

Finally,

\begin{align}
M_1(r)=M^{1}(r)+M^{2}(r)+M^{3}(r)=0.
\end{align}
 \begin{flushright}$\blacksquare$\end{flushright}

 In the following we consider a much more general
class of discontinuous functions
\begin{equation}\label{nr1g}
g(x,y,\varepsilon)=\alpha_{i}h'(x)+\varepsilon y \ \ if \
x\in(a_{i}, a_{i+1}], \ i=0,1,...,n
\end{equation}
\noindent with the same convention for the last interval, i.e. it
is $(a_n,a_{n+1}),$ where $h'(x)$ is the derivative of a
differentiable function $h(x)$ and we will prove a similar result
given by:
\begin{prop}\label{p2}
If $g$ is a non-smooth function given by (\ref{nr1g}), then, for
any $0<\varepsilon\ll 1$ and any two sequences
$(a_{i})_{i=0,1,...,n+1}$ and $(\alpha_{i})_{i=0,1,...,n}$ as
above, the  first Melnikov function of the system (\ref{nr1}) is
identically zero.
\end{prop}
The proof is similar. Proceeding as above, we have that:
\begin{align}
M^{1}(r):&= \int_{0}^{\pi/2}r\cos t\ g(r\sin t)dt
=\sum^{m}_{i=0}\int^{t_{i+1}}_{t_{i}}r\alpha_{i}\cos
t\ h'(r\sin t)dt\\
&=h(r)\alpha_{m}-h(0)\alpha_{0}-\sum^{m}_{i=1}(\alpha_{i}-\alpha_{i-1})h(a_{i})\nonumber,
\end{align}
\begin{align}
M^{2}(r):&= \int_{0}^{\pi/2}r\cos t\ g(r\sin t)dt
=\sum^{m}_{i=0}\int^{t_{i}}_{t_{i+1}}r\alpha_{i}\cos
t\ h'(r\sin t)dt\\
&=-h(r)\alpha_{m}+h(-r)\alpha_{0}+\sum^{m}_{i=1}(\alpha_{i}-\alpha_{i-1})h(a_{i})\nonumber,
\end{align}
and
\begin{align}
M^{3}(r):= \int_{3\pi/2}^{2\pi}r\cos t\ g(r\sin t)dt=
\int_{3\pi/2}^{2\pi}r\alpha_0\cos t h'(r\sin
t)dt=h(0)\alpha_{0}-h(-r)\alpha_{0},
\end{align}
so we arrive to the same result:
\begin{align}
M_1(r)=M^{1}(r)+M^{2}(r)+M^{3}(r)=0.
\end{align}
\begin{flushright}$\blacksquare$\end{flushright}

As the first Melnikov function is identically zero, we can say
nothing about the number of limit cycles. However, we conjecture
the following fact:\\
\\
CONJECTURE: \emph{The non-smooth system (\ref{nr1}) with the
function $g$ given by (\ref{nr1g}) has no limit cycles.}\\
\\
We base this conjecture of the following fact. Compute the second
Melnikov function for the system (\ref{nr1}) with $g$ given by
(\ref{nr1g}). One can observe that $\frac{\partial f}{\partial
x}(x,y,0)+\frac{\partial g}{\partial y}(x,y,0)=0$ on any strip
where $g$ is defined, so using (\ref{nr05}) we get that:
\begin{equation}\label{nr5}
M_2(r)=\oint_{C_r}\frac{\partial g}{\partial
\varepsilon}(x,y,0)dx-\frac{\partial f}{\partial
\varepsilon}(x,y,0)dy=\oint_{C_r}ydx=\int_{0}^{2\pi}r^2\cos^2
tdt=r^2\pi\neq 0.
\end{equation}
As $M_2(r)=0,r>0$ has no root the above Conjecture is justified by
Theorem \ref{teo1}. However, we can not present this Conjecture as
a result, because the formula of the second Melnikov function
which we used from \cite{ile1} was proved for continuous functions
but we applied it for functions with a finite number of
discontinuities.

\section{Conclusions}

In this work we started to explore an example of a class of
non-smooth dynamical systems with the first Melnikov function
identically zero. We applied the formula of the second Melnikov
function reported in a work of Iliev \cite{ile1} in order to
investigate further the existence of limit cycles. This formula
has been deduced for continuous functions but it is quite
plausible it remains valid for piecewise functions.

\end{document}